\documentclass[12pt]{amsart}

\usepackage{enumitem}
\usepackage{psfrag}
\usepackage{graphicx}
\usepackage{pinlabel}
\usepackage{graphicx}
\usepackage{amsmath}
\usepackage{amssymb}
\usepackage{amscd}
\usepackage{times}
\usepackage{pb-diagram}
\usepackage{graphicx}
\usepackage{wrapfig}
\usepackage{xspace}
\usepackage{hyperref}
\usepackage{url}
\usepackage{caption}
\usepackage{subcaption}
\usepackage{fancyhdr}
\usepackage{overpic}
\usepackage{color}
\usepackage[square,comma,sort&compress,numbers]{natbib}
\usepackage{latexsym}
\usepackage{mathrsfs}
\usepackage{amsfonts}
\usepackage{hyperref}
\usepackage{amsthm}
\usepackage{appendix}
\usepackage{pdfsync}

\theoremstyle{definition}
\newtheorem{theorem}{Theorem}[section]
\newtheorem{lemma}[theorem]{Lemma}
\newtheorem{proposition}[theorem]{Proposition}

\newtheorem{definition}[theorem]{Definition}

\newtheorem{remark}[theorem]{Remark}
\newtheorem*{theorem*}{Theorem}
\title[Spherical Combinatorial Curvature Flows on Infinitely Triangulated Surfaces]
{Combinatorial Curvature Flows for spherical vertex scalings on Infinitely Triangulated Surfaces}
\author{Yuerong Bian, Xu Xu}
\date{\today}

\address{School of Mathematics and Statistics, Wuhan University, Wuhan, 430072, P. R. China}
\email{yuerongbian@whu.edu.cn}

\address{School of Mathematics and Statistics, Wuhan University, Wuhan, 430072, P. R. China}
\email{xuxu2@whu.edu.cn}



\begin{document}

\begin{abstract}
We study the combinatorial Yamabe and Calabi flows for spherical vertex scalings on infinitely triangulated surfaces. 
Under the assumptions of uniformly bounded vertex degree and an $\epsilon$-uniformly nondegenerate initial
metric, we establish the short-time existence of smooth solutions to both combinatorial curvature flows.
We also establish short-time uniqueness for the combinatorial Yamabe flow.
For the combinatorial Calabi flow, we establish short-time uniqueness among solutions that are continuous in time with respect to the $L^\infty(V)$ norm.
\end{abstract}
	\maketitle
	
	\section{Introduction}
Finding metrics with prescribed curvature is a basic problem in geometry. 
Curvature flows provide a method for studying this problem, with Hamilton's Ricci flow \cite{Hamilton1982} and Calabi's flow \cite{Calabi1982,Calabi1985} as important examples.
Combinatorial Ricci flow (Yamabe flow) and combinatorial Calabi flow are discrete analogues of Hamilton's Ricci flow and the Calabi flow.
Chow and Luo \cite{ChowLuo2003} introduced the combinatorial Ricci flow for Thurston's circle packing
metrics \cite{Thurston1980} on finitely triangulated closed surfaces. Luo \cite{Luo2004} introduced the combinatorial Yamabe flow for Euclidean vertex scalings on finitely triangulated closed surfaces. 
Ge \cite{GeThesis2012} then introduced the combinatorial Calabi flow for Thurston's Euclidean circle packing metrics and vertex scalings on finitely triangulated closed surfaces. 
There has been extensive research on combinatorial Ricci flow (Yamabe flow) and combinatorial Calabi flow on closed surfaces in Euclidean and hyperbolic background geometry. See \cite{GeXu2016,GeHua2018, ZhuXu2019, BobenkoPinkallSpringborn2015, GuLuoSunWu2018, GuGuoLuoSunWu2018} and others.
Research on spherical vertex scalings is relatively rare in the literature. 
Izmestiev, Prosanov and Wu
\cite{IzmestievProsanovWu2024} recently studied the prescribed positive discrete curvature
problem for spherical vertex scalings on the $2$-sphere.
Nie \cite{Nie} studied the spherical circle packings using combinatorial geodesic curvature. 
Ge, Hua and Zhou \cite{GeHuaZhouInfinite2024JFA} studied the corresponding combinatorial Ricci flow for spherical circle packings using the combinatorial geodesic curvatures.
Izmestiev and Lam \cite{IzmestievLam2025} studied the spherical discrete Laplacians on triangulated surfaces closely related to spherical vertex scalings.
	
Recently, several works studied combinatorial curvature flows on infinitely triangulated surfaces in Euclidean and hyperbolic background geometry. Ge, Hua and Zhou
	\cite{GeHuaZhouInfinite2025} studied the combinatorial Ricci flows for Thurston's circle packing metrics on infinitely triangulated surfaces in Euclidean and hyperbolic background geometries. 
	For vertex scalings, Ji \cite{Ji2025} studied the Euclidean combinatorial Yamabe flow on infinitely triangulated surfaces. 
In \cite{BianXuYamabe2026}, we established the short-time existence and uniqueness for the corresponding hyperbolic Yamabe flow under the assumptions of uniformly bounded vertex degree and
an $\epsilon$-uniformly nondegenerate initial metric.
	
In this paper, we study the combinatorial Yamabe and Calabi flows for spherical vertex scalings on infinitely triangulated surfaces.
If the vertex degrees are bounded by $D$ and the initial
piecewise spherical metric is $\epsilon$-uniformly nondegenerate, we prove that both combinatorial curvature flows admit solutions that are smooth in time on intervals depending only on $\epsilon$ and $D$. 
We also prove the short-time uniqueness. 
For the combinatorial Calabi flow, uniqueness holds among solutions that are continuous in time in the $L^\infty$
norm; the solutions obtained in the existence proof have this property.
	
The proofs use uniform estimates for the curvature derivatives and their variations under spherical vertex scalings. 
An angle estimate ensures that the metrics remain uniformly nondegenerate when the conformal factors are sufficiently small in the $L^\infty$ norm.
These estimates give a common existence interval
for solutions on finite subcomplexes. 
We then obtain solutions on the infinitely triangulated surface by Arzel\`a--Ascoli theorem and
a diagonal argument. 
Uniqueness follows from uniform Lipschitz estimates and Gr\"onwall's inequality.

Throughout this paper, $(S,\mathcal T)$ denotes a connected,
locally finite, infinitely triangulated surface without
boundary, where $\mathcal T=(V,E,F)$ is the triangulation of the surface with $V, E, F$ being the sets of vertices, edges and faces respectively.	
We use $i$, $\{ij\}$, and
$\triangle ijk$ to denote a vertex, an edge, and a face,
respectively. We write $j\sim i$ if $\{ij\}\in E$ and denote
the degree of $i$ by $\deg(i)$.
A piecewise spherical (PS) metric on $\mathcal T$ is determined
by a length function $d\colon E\to(0,\pi)$ such that, for each
face $\triangle ijk\in F$,
\[
d_{ij}+d_{ik}>d_{jk}\quad\text{(cyclically)},
\qquad
d_{ij}+d_{ik}+d_{jk}<2\pi.
\]
These lengths determine, up to isometry, nondegenerate
spherical triangles on the unit sphere. Gluing these triangles
isometrically along their common edges gives a spherical cone
metric on $S$. We identify the spherical cone metric with the edge length
function $d$.
The combinatorial curvature at $i\in V$ is defined by
\[
K_i=2\pi-\sum_{\triangle ijk\in F}\theta_i^{jk},
\]
where $\theta_i^{jk}$ is the interior angle at $i$ of the triangle $\triangle ijk$.

Motivated by Luo's definition of vertex scalings of piecewise linear metrics \cite{Luo2004}, 
Bobenko-Pinkall-Springborn \cite{BobenkoPinkallSpringborn2015}
introduced the following definition of vertex scalings of PS metrics. 

\begin{definition}[\cite{BobenkoPinkallSpringborn2015}]\label{Def: vertex_scaling}
	Let $(S,\mathcal T)$ be an infinitely triangulated surface,
	and let $d$ and $\tilde d$ be two PS metrics on $\mathcal T$.
	We say that $\tilde d$ is a spherical vertex scaling of $d$ if
	there exists a function $u\colon V\to\mathbb R$ such that
	\begin{equation}\label{eq: vertex_scaling}
		\sin\frac{\tilde d_{ij}}{2}
		=e^{\frac{u_i+u_j}{2}}\sin\frac{d_{ij}}{2},
		\qquad \{ij\}\in E.
	\end{equation}
	The function $u$ is called a discrete conformal factor, and we
	write $\tilde d=u*d$.
\end{definition}

The reader may refer to \cite{BobenkoPinkallSpringborn2015} for a nice explanation of Definition \ref{Def: vertex_scaling} in terms of 3-dimensional hyperbolic geometry. The combinatorial Yamabe flow and combinatorial Calabi flow for spherical vertex scalings are defined as follows.

\begin{definition}\label{Def: CYF}
	Let $(S,\mathcal T)$ be an infinitely triangulated surface
	with an initial PS metric $d_0$.
	The spherical combinatorial Yamabe flow is defined by
\begin{equation}\label{Eq: CYF}
	\begin{cases}
		\dfrac{du_i}{dt}=-K_i,\\[4pt]
		u_i(0)=0,
	\end{cases}
	\qquad i\in V,
\end{equation}
	where $K_i$ is the combinatorial curvature of the PS metric
	$u(t)*d_0$ at $i$.
\end{definition}

\begin{definition}\label{Def:VCLBF}
	Let $(S,\mathcal T)$ be an infinitely triangulated surface
	with an initial PS metric $d_0$.
	The spherical combinatorial Calabi flow is defined by
\begin{equation}\label{Eq:VCLBF}
	\begin{cases}
		\dfrac{du_i}{dt}=(\Delta K)_i,\\[4pt]
		u_i(0)=0,
	\end{cases}
	\qquad i\in V,
\end{equation}
	where the discrete Laplace operator is given by
	\[
	(\Delta f)_i
	=-\sum_{j\in V}\frac{\partial K_i}{\partial u_j}f_j,
	\qquad f\colon V\to\mathbb R.
	\]
Here $K_i$ is the combinatorial curvature of the PS metric $u(t)*d_0$ at  $i$, and the sum in the discrete Laplace operator is finite for each \(i\), since
\(\partial K_i/\partial u_j=0\) unless \(j=i\) or \(j\sim i\).
\end{definition}

To state our main results, we introduce an
$\epsilon$-uniform nondegeneracy condition for PS metrics.
This condition is motivated by the description of the space
of spherical triangles in terms of their interior angles
\cite[Section~2.1]{Luo2006}.
For $(\theta_1,\theta_2,\theta_3)\in(0,\pi)^3$, set
\[
\theta_i^*:=\frac12(\pi+\theta_i-\theta_j-\theta_k),
\qquad \{i,j,k\}=\{1,2,3\},
\]
and let
\[
M:=\left\{
(\theta_1,\theta_2,\theta_3)\in(0,\pi)^3:
\theta_i^*>0\ (i=1,2,3),\quad
\theta_1+\theta_2+\theta_3>\pi
\right\}.
\]
The set $M$ consists of the angle triples of nondegenerate
spherical triangles on the unit sphere.
For $\epsilon>0$, define the compact subset
\[
\mathcal K_\epsilon^S:=
\left\{
(\theta_1,\theta_2,\theta_3)\in M:
\begin{array}{l}
	\epsilon\leq\theta_i\leq\pi-\epsilon
	\quad (i=1,2,3),\\
	\theta_1+\theta_2+\theta_3\geq\pi+\epsilon,\\
	\theta_i+\theta_j-\theta_k\leq\pi-\epsilon
	\quad (\{i,j,k\}=\{1,2,3\})
\end{array}
\right\}.
\]

\begin{definition}\label{Def: CRF}
	Let $(S,\mathcal T)$ be an infinitely triangulated surface
	with a PS metric $d$.
	For $\epsilon>0$, we say that $d$ is
	$\epsilon$-uniformly nondegenerate if
	\[
	\bigl(\theta_i^{jk},\theta_j^{ik},\theta_k^{ij}\bigr)
	\in\mathcal K_\epsilon^S
	\]
	for every face $\triangle ijk\in F$.
\end{definition}

We denote by $L^\infty(V)$ the space of bounded functions on $V$,
with norm
\[
\|u\|_{L^\infty(V)}=\sup_{i\in V}|u_i|.
\]
The notation $u\in C_t^\infty(V\times[0,T])$ means that
$u_i(t)\in C^\infty([0,T])$ for every $i\in V$.
Continuity in $C([0,T];L^\infty(V))$ is understood with respect
to this norm.

In this paper, we establish the following short-time existence result for the spherical combinatorial Yamabe flow and Calabi flow on infinitely triangulated surfaces.

\begin{theorem}[Short-time existence]\label{thm:vxste}
	Let $(S,\mathcal T)$ be an infinitely triangulated surface with
	vertex degree bounded by $D$, and let $d_0$ be an
	$\epsilon$-uniformly nondegenerate initial PS metric on
	$\mathcal T$.
	\begin{enumerate}
		\item[(a)] There exists $T_1=T_1(\epsilon,D)>0$ such that
		the spherical combinatorial Yamabe flow \eqref{Eq: CYF}
		admits a solution $u\in C_t^\infty(V\times[0,T_1])$.
		
		\item[(b)] There exists $T_2=T_2(\epsilon,D)>0$ such that
		the spherical combinatorial Calabi flow \eqref{Eq:VCLBF}
		admits a solution $u\in C_t^\infty(V\times[0,T_2])$.
	\end{enumerate}
\end{theorem}

The solution constructed in Theorem~\ref{thm:vxste}(b) is also
continuous in the $L^\infty$ norm and hence satisfies the
continuity condition in Theorem~\ref{thm:uq}(b) below.
Under the same assumptions, we further establish the following short-time
uniqueness for both combinatorial curvature flows.

\begin{theorem}[Short-time uniqueness]\label{thm:uq}
	Let $(S,\mathcal T)$ be an infinitely triangulated surface with
	vertex degree bounded by $D$, and let $d_0$ be an
	$\epsilon$-uniformly nondegenerate initial PS metric on
	$\mathcal T$.
	Let $T_1$ and $T_2$ be the existence times given by
	Theorem~\ref{thm:vxste}.
	\begin{enumerate}
		\item[(a)] Let $u,v\in C_t^\infty(V\times[0,T])$ be two
		solutions of the spherical combinatorial Yamabe flow
		\eqref{Eq: CYF}, where $0<T\leq T_1$.
		Then $u(t)=v(t)$ for all $t\in[0,T]$.
		
		\item[(b)] Let
		\[
		u,v\in C_t^\infty(V\times[0,T])
		\cap C([0,T];L^\infty(V))
		\]
		be two solutions of the spherical combinatorial Calabi
		flow \eqref{Eq:VCLBF}, where $0<T\leq T_2$.
		Then $u(t)=v(t)$ for all $t\in[0,T]$.
	\end{enumerate}
\end{theorem}

The paper is organized as follows. In Section~2, we derive
the variational formulas for spherical vertex scalings and
establish the uniform estimates used in the proofs.
In Sections~3 and~4, we prove the short-time existence and
uniqueness for the spherical combinatorial Yamabe and Calabi
flows, respectively.
	\\
	\\
	\textbf{Acknowledgments:} 
	The research of Xu Xu is supported by the National Natural Science Foundation of China under Grant No. 12471057 and Hubei Provincial Natural Science Fund for Outstanding Research
	Groups under Grant No. 2026AFA008. 
	\\
	\\
	\noindent\textbf{AI Statement:} All mathematical statements and proofs in this paper were independently developed by the authors, who take full responsibility for the content.
	AI tools were only used for English-language editing.

\section{Preliminaries}
\subsection{Variational formulas for spherical vertex scalings}
We first derive the angle variational formula for spherical
vertex scalings and then give the corresponding variational
formulas for the combinatorial curvatures.
We need the following formula to derive the variational formulas.

\begin{lemma}\label{lem:spherical-half-angle-identity}
	Let $d_{jk},d_{ik},d_{ij}$ be the edge lengths of a nondegenerate
	spherical triangle $\triangle ijk$, and let
	$\theta_i^{jk},\theta_j^{ik},\theta_k^{ij}$ be the corresponding
	interior angles at the vertices $i, j, k$ respectively. Then
	\[
	\frac{
		\sin^2\frac{d_{jk}}{2}
		+\sin^2\frac{d_{ik}}{2}
		-\sin^2\frac{d_{ij}}{2}
	}{
		\sin\frac{d_{jk}}{2}\sin\frac{d_{ik}}{2}
	}
	=
	2\sin\frac{\theta_i^{jk}+\theta_j^{ik}-\theta_k^{ij}}{2}
	\cos\frac{d_{ij}}{2}.
	\]
\end{lemma}
\begin{proof}
	We adapt the proof of \cite[Lemma~11]{GuGuoLuoSunWu2018}
	to the spherical case.
	By the spherical cosine law in terms of the angles,
	\[
	\cos d_{jk}
	=
	\frac{
		\cos\theta_i^{jk}
		+\cos\theta_j^{ik}\cos\theta_k^{ij}
	}{
		\sin\theta_j^{ik}\sin\theta_k^{ij}
	},
	\]
	and the analogous formulas hold for $\cos d_{ik}$ and $\cos d_{ij}$.
	Therefore,
	\begin{equation}\label{eq:spherical-half-angle 0}
	\begin{aligned}
		& 	\sin^2\frac{d_{jk}}{2}+\sin^2\frac{d_{ik}}{2}-\sin^2\frac{d_{ij}}{2} \\
		= &\frac{1}{2} (1-\cos d_{jk} -\cos d_{ik} + \cos d_{ij}  ) \\
		=& \frac{1}{2} \left[1-
		\frac{\cos \theta_i^{jk} + \cos \theta_j^{ik}\cos \theta_k^{ij}}{\sin \theta_j^{ik}\sin \theta_k^{ij}}
		-\frac{\cos \theta_j^{ik} + \cos \theta_i^{jk}\cos \theta_k^{ij}}{\sin \theta_i^{jk}\sin \theta_k^{ij}}
		+\frac{\cos \theta_k^{ij} + \cos \theta_i^{jk}\cos \theta_j^{ik}}{\sin \theta_i^{jk}\sin \theta_j^{ik}}
		\right] \\
		=&- \frac{1}{2\sin \theta_i^{jk}\sin \theta_j^{ik}\sin \theta_k^{ij}}
		\big(\sin (\theta_i^{jk} + \theta_j^{ik}) - \sin \theta_k^{ij}\big)\big(\cos \theta_k^{ij} + \cos (\theta_i^{jk} - \theta_j^{ik})\big) \\
		= &-\frac{2\sin \frac{\theta_i^{jk} + \theta_j^{ik} - \theta_k^{ij}}{2}
			\cos \frac{\theta_i^{jk} + \theta_j^{ik} + \theta_k^{ij}}{2}
			\cos \frac{\theta_i^{jk} + \theta_k^{ij} - \theta_j^{ik}}{2}
			\cos \frac{\theta_j^{ik} + \theta_k^{ij} - \theta_i^{jk}}{2}}
		{\sin \theta_i^{jk}\sin \theta_j^{ik}\sin \theta_k^{ij}}.
	\end{aligned}
\end{equation}
	On the other hand, the half-angle formulas give
	
	\begin{equation}\label{eq:spherical-half-angle}
		\begin{aligned}
			\sin^2\frac{d_{jk}}2
			&= \frac{1}{2}(1-\cos d_{jk}) \\
			&= \frac{1}{2}\left(1-
			\frac{\cos \theta_i^{jk} + \cos \theta_j^{ik}\cos \theta_k^{ij}}
			{\sin \theta_j^{ik}\sin \theta_k^{ij}}
			\right) \\
			&= -\frac{1}{2}
			\frac{\cos \theta_i^{jk} + \cos(\theta_j^{ik} + \theta_k^{ij})}
			{\sin \theta_j^{ik}\sin \theta_k^{ij}} \\
			&= -\frac{
				\cos \frac{\theta_i^{jk}+\theta_j^{ik}+\theta_k^{ij}}{2}
				\cos \frac{\theta_i^{jk}-\theta_j^{ik}-\theta_k^{ij}}{2}
			}{
				\sin \theta_j^{ik}\sin \theta_k^{ij}
			},
		\end{aligned}
	\end{equation}
	and
	\begin{equation}\label{eq:spherical-half-angle 2}
		\begin{aligned}
	\sin^2\frac{d_{ik}}2
	=
	-\frac{
		\cos\frac{\theta_i^{jk}+\theta_j^{ik}+\theta_k^{ij}}2
		\cos\frac{\theta_j^{ik}-\theta_i^{jk}-\theta_k^{ij}}2
	}{
		\sin\theta_i^{jk}\sin\theta_k^{ij}
	}.
		\end{aligned}
	\end{equation}
	Since
$\theta_i^{jk}+\theta_j^{ik}+\theta_k^{ij}>\pi$,
	we have
	$\cos\frac{\theta_i^{jk}+\theta_j^{ik}+\theta_k^{ij}}2<0$.
	The spherical cosine law also gives
	
	\begin{equation}\label{eq:spherical-cos-half-angle}
		\begin{aligned}
			\cos^2 \frac{d_{ij}}{2}
			&= \frac{1}{2}\left(
			\frac{\cos \theta_k^{ij} + \cos \theta_i^{jk}\cos \theta_j^{ik}}
			{\sin \theta_i^{jk}\sin \theta_j^{ik}}
			+1
			\right) \\
			&= \frac{1}{2}
			\frac{\cos \theta_k^{ij} + \cos(\theta_i^{jk} - \theta_j^{ik})}
			{\sin \theta_i^{jk}\sin \theta_j^{ik}} \\
			&=
			\frac{
				\cos \frac{\theta_i^{jk}+\theta_k^{ij}-\theta_j^{ik}}{2}
				\cos \frac{\theta_j^{ik}+\theta_k^{ij}-\theta_i^{jk}}{2}
			}{
				\sin \theta_i^{jk}\sin \theta_j^{ik}
			}.
		\end{aligned}
	\end{equation}
	Since $0<d_{ij}<\pi$, we have $\cos(d_{ij}/2)>0$.
	Taking positive square roots, we obtain
	\[
	\begin{aligned}
		& \frac{
			\sin^2(d_{jk}/2) + \sin^2(d_{ik}/2) - \sin^2(d_{ij}/2)
		}{
			\sin(d_{jk}/2)\sin(d_{ik}/2)
		} \\
		&= 
		-\frac{
			2\sin \frac{\theta_i^{jk}+\theta_j^{ik}-\theta_k^{ij}}{2}
			\cos \frac{\theta_i^{jk}+\theta_j^{ik}+\theta_k^{ij}}{2}
			\cos \frac{\theta_i^{jk}+\theta_k^{ij}-\theta_j^{ik}}{2}
			\cos \frac{\theta_j^{ik}+\theta_k^{ij}-\theta_i^{jk}}{2}
		}{
			\sin \theta_i^{jk}\sin \theta_j^{ik}\sin \theta_k^{ij}
			\cdot
			\sqrt{
				-\frac{
					\cos \frac{\theta_i^{jk}+\theta_j^{ik}+\theta_k^{ij}}{2}
					\cos \frac{\theta_i^{jk}-\theta_j^{ik}-\theta_k^{ij}}{2}
				}{
					\sin \theta_j^{ik}\sin \theta_k^{ij}
				}
			}
			\sqrt{
				-\frac{
					\cos \frac{\theta_i^{jk}+\theta_j^{ik}+\theta_k^{ij}}{2}
					\cos \frac{\theta_j^{ik}-\theta_i^{jk}-\theta_k^{ij}}{2}
				}{
					\sin \theta_i^{jk}\sin \theta_k^{ij}
				}
			}
		} \\
		&= 2\sin \frac{\theta_i^{jk} + \theta_j^{ik} - \theta_k^{ij}}{2}
		\cdot
		\sqrt{
			\frac{
				\cos \frac{\theta_i^{jk}+\theta_k^{ij}-\theta_j^{ik}}{2}
				\cos \frac{\theta_j^{ik}+\theta_k^{ij}-\theta_i^{jk}}{2}
			}{
				\sin \theta_i^{jk}\sin \theta_j^{ik}
			}
		} \\
		&= 2\sin \frac{\theta_i^{jk} + \theta_j^{ik} - \theta_k^{ij}}{2}
		\cdot \cos \frac{d_{ij}}{2}.
	\end{aligned}
	\]
where \eqref{eq:spherical-half-angle 0}, \eqref{eq:spherical-half-angle}, \eqref{eq:spherical-half-angle 2} are used in the second equality, and \eqref{eq:spherical-cos-half-angle} is used in the last equality.
\end{proof}

\begin{lemma}\label{lem:spherical-angle-variation}
	Let $d_0$ be a PS metric on $\triangle ijk$, and let $u$ be a
	discrete conformal factor such that $d=u*d_0$ defines a
	nondegenerate spherical triangle. Denote its edge lengths by
	$d_{jk},d_{ik},d_{ij}$ and the corresponding opposite angles by
	$\theta_i^{jk},\theta_j^{ik},\theta_k^{ij}$. Then
	\[
	\frac{\partial\theta_i^{jk}}{\partial u_j}
	=\frac{\partial\theta_j^{ik}}{\partial u_i}
	=\frac{1}{2\cos^2\frac{d_{ij}}{2}}
	\tan\frac{\theta_i^{jk}+\theta_j^{ik}-\theta_k^{ij}}2.
	\]
\end{lemma}

\begin{proof}
	We adapt the argument of Xu and Zheng
	\cite[Lemma~A.1]{XuZheng2022} to the spherical case.
	By the derivative cosine law
	(see \cite[Lemma~A.1]{ChowLuo2003}), we have
	\[
	\frac{\partial\theta_i^{jk}}{\partial d_{jk}}
	=\frac{\sin d_{jk}}{A},
	\qquad
	\frac{\partial\theta_i^{jk}}{\partial d_{ik}}
	=-\frac{\sin d_{jk}\cos\theta_k^{ij}}{A},
	\qquad
	\frac{\partial\theta_i^{jk}}{\partial d_{ij}}
	=-\frac{\sin d_{jk}\cos\theta_j^{ik}}{A},
	\]
	where $A=\sin d_{ik}\sin d_{ij}\sin\theta_i^{jk}$.
	By the vertex-scaling relationship \eqref{eq: vertex_scaling},
	\[
	\frac{\partial d_{jk}}{\partial u_i}=0,
	\qquad
	\frac{\partial d_{jk}}{\partial u_j}
	=\frac{\partial d_{jk}}{\partial u_k}
	=\tan\frac{d_{jk}}2.
	\]
	Thus, by the chain rule and the spherical cosine law, we have
	
	\begin{equation}\label{eq:spherical-angle-variation-lengths}
		\begin{aligned}
			\frac{\partial\theta_i^{jk}}{\partial u_j}
			&=
			\frac{\sin d_{jk}}{A}
			\left(
			\tan\frac{d_{jk}}2
			-\cos\theta_j^{ik}\tan\frac{d_{ij}}2
			\right)\\
			&=
			\frac{1-\cos d_{jk}}{A}
			-
			\frac{\cos d_{ik}-\cos d_{ij}\cos d_{jk}}
			{A(1+\cos d_{ij})}\\
			&=
			\frac{1+\cos d_{ij}-\cos d_{jk}-\cos d_{ik}}
			{A(1+\cos d_{ij})}.
		\end{aligned}
	\end{equation}
Applying \eqref{eq:spherical-cos-half-angle} with
$(i,j,k)$ replaced by $(j,k,i)$ and $(k,i,j)$ gives
the corresponding formulas for $d_{jk}$ and $d_{ik}$.
Since
$
-\pi<\theta_i^{jk}+\theta_j^{ik}-\theta_k^{ij}<\pi,
$
we have
$\cos\frac{\theta_i^{jk}+\theta_j^{ik}-\theta_k^{ij}}2>0$.
Taking positive square roots and using
\eqref{eq:spherical-cos-half-angle} again, we obtain
\begin{equation}\label{eq:spherical-half-angle-product}
    \begin{aligned}
        &\cos\frac{d_{jk}}2\cos\frac{d_{ik}}2
        \sin\theta_k^{ij}\\
        &\quad=
        \cos\frac{\theta_i^{jk}+\theta_j^{ik}-\theta_k^{ij}}2
        \sqrt{
            \frac{
                \cos\frac{\theta_i^{jk}-\theta_j^{ik}
                    +\theta_k^{ij}}2
                \cos\frac{-\theta_i^{jk}+\theta_j^{ik}
                    +\theta_k^{ij}}2
            }{
                \sin\theta_i^{jk}\sin\theta_j^{ik}
            }
        }\\
        &\quad=
        \cos\frac{d_{ij}}2
        \cos\frac{\theta_i^{jk}+\theta_j^{ik}-\theta_k^{ij}}2.
    \end{aligned}
\end{equation}
	By the spherical sine law, we have 
	$A=\sin d_{ik}\sin d_{jk}\sin\theta_k^{ij}$.
	Hence \eqref{eq:spherical-angle-variation-lengths} yields
	\begin{equation}\label{Eq: derivarive of theta_i of u_j}
	\begin{aligned}
		\frac{\partial\theta_i^{jk}}{\partial u_j}
		&=
		\frac{
			1+\cos d_{ij}-\cos d_{jk}-\cos d_{ik}
		}{
			\sin d_{ik}\sin d_{jk}\sin\theta_k^{ij}
			(1+\cos d_{ij})
		} \\
		&=
		\frac{
			2\left(
			\sin^2\frac{d_{jk}}{2}
			+\sin^2\frac{d_{ik}}{2}
			-\sin^2\frac{d_{ij}}{2}
			\right)
		}{
			\sin d_{ik}\sin d_{jk}\sin\theta_k^{ij}
			(1+\cos d_{ij})
		} \\
		&=
		\frac{
			4\sin\frac{d_{jk}}{2}\sin\frac{d_{ik}}{2}
			\cos\frac{d_{ij}}{2}
		}{
			\sin d_{ik}\sin d_{jk}\sin\theta_k^{ij}
			(1+\cos d_{ij})
		}
		\sin\frac{\theta_i^{jk}+\theta_j^{ik}-\theta_k^{ij}}{2} \\
		&=
		\frac{
			\sin\frac{\theta_i^{jk}+\theta_j^{ik}-\theta_k^{ij}}{2}
		}{
			2\cos\frac{d_{ij}}{2}
			\cos\frac{d_{jk}}{2}
			\cos\frac{d_{ik}}{2}
			\sin\theta_k^{ij}
		} \\
		&=
		\frac{1}{2\cos^2\frac{d_{ij}}{2}}
		\tan\frac{\theta_i^{jk}+\theta_j^{ik}-\theta_k^{ij}}{2}.
	\end{aligned}
	\end{equation}
	The third equality follows from
	Lemma~\ref{lem:spherical-half-angle-identity}, and the last equality
	follows from \eqref{eq:spherical-half-angle-product}.
  The symmetry follows by interchanging $i$ and $j$ in \eqref{Eq: derivarive of theta_i of u_j}.
\end{proof}

Let $d=u*d_0$ be a PS metric on $\mathcal T$.
For each edge $\{ij\}\in E$, let $\triangle ijk_1$ and
$\triangle ijk_2$ be two faces  incident to the edge $\{ij\}$. Set
\begin{equation}\label{eq:spherical-coefficients-definition}
    \begin{aligned}
        B_{ij}
        &:=\frac{\partial\theta_i^{jk_1}}{\partial u_j}
        +\frac{\partial\theta_i^{jk_2}}{\partial u_j},\\
        A_i
        &:=\sum_{\triangle ijk\in F}
        \frac{\partial\operatorname{Area}(\triangle ijk)}{\partial u_i}.
    \end{aligned}
\end{equation}

\begin{proposition}
	\label{prop:spherical-curvature-variation}
	For $j\sim i$, we have
	\begin{equation}\label{eq:spherical-curvature-off-diagonal}
		\begin{aligned}
			\frac{\partial K_i}{\partial u_j}
			&=\frac{\partial K_j}{\partial u_i}=-B_{ij}=
			-\frac{1}{2\cos^2\frac{d_{ij}}{2}}
			\sum_{\ell=1}^{2}
			\tan\frac{
				\theta_i^{jk_\ell}+\theta_j^{ik_\ell}-\theta_{k_\ell}^{ij}
			}{2}.
		\end{aligned}
	\end{equation}
	Moreover,
	\begin{equation}\label{eq:spherical-curvature-A-B}
		\begin{aligned}
			A_i
			&=\sum_{j\sim i}B_{ij}(1-\cos d_{ij})=\sum_{j\sim i}\tan^2\frac{d_{ij}}{2}
			\sum_{\ell=1}^{2}
			\tan\frac{
				\theta_i^{jk_\ell}+\theta_j^{ik_\ell}-\theta_{k_\ell}^{ij}
			}{2},
		\end{aligned}
	\end{equation}
	and
	
\begin{equation}\label{eq:spherical-curvature-diagonal}
		\begin{aligned}
			\frac{\partial K_i}{\partial u_i}
			&=-A_i+\sum_{j\sim i}B_{ij}
			=\sum_{j\sim i}B_{ij}\cos d_{ij}=\sum_{j\sim i}
			\frac{\cos d_{ij}}{2\cos^2\frac{d_{ij}}{2}}
			\sum_{\ell=1}^{2}
			\tan\frac{
				\theta_i^{jk_\ell}+\theta_j^{ik_\ell}-\theta_{k_\ell}^{ij}
			}{2}.
		\end{aligned}
	\end{equation}
	If $j\ne i$ and $j\not\sim i$, then
	$\partial K_i/\partial u_j=0$.
\end{proposition}

\begin{proof}
	Formula \eqref{eq:spherical-curvature-off-diagonal} and the
	vanishing of $\partial K_i/\partial u_j$ for $j\ne i$ and
	$j\not\sim i$ follow from the definition of $K_i$ and
	Lemma~\ref{lem:spherical-angle-variation}.
Using the identity
$\operatorname{Area}(\triangle ijk)
=\theta_i^{jk}+\theta_j^{ik}+\theta_k^{ij}-\pi$
and the symmetry of the angle derivatives, we obtain
\begin{equation}\label{eq: derivative of K w.r.t. u}
	\begin{aligned}
    		\frac{\partial K_i}{\partial u_i}
		&=-\sum_{\triangle ijk\in F}
		\frac{\partial\theta_i^{jk}}{\partial u_i}=\sum_{\triangle ijk\in F}
		\left(
		\frac{\partial\theta_i^{jk}}{\partial u_j}
		+
		\frac{\partial\theta_i^{jk}}{\partial u_k}
		\right)-A_i=\sum_{j\sim i}B_{ij}-A_i.
	\end{aligned}
\end{equation}

	By the spherical area variational formula of Glickenstein and Thomas \cite[Proposition 11]{GlickensteinThomas2017},
	\begin{equation}\label{eq:ddarea}
		\frac{\partial\operatorname{Area}(\triangle ijk)}{\partial u_i}
		=
		(1-\cos d_{ij})\frac{\partial\theta_j^{ik}}{\partial u_i}
		+
		(1-\cos d_{ik})\frac{\partial\theta_k^{ij}}{\partial u_i}.
	\end{equation}
	Summing \eqref{eq:ddarea} over the faces incident to $i$, using
	the symmetry of the angle derivatives in 	Lemma~\ref{lem:spherical-angle-variation}, and grouping terms by
	the edges incident to $i$, we obtain
	\begin{equation}\label{eq: Ai in Bij}
	\begin{aligned}
		A_i
		&=\sum_{\triangle ijk\in F}
		\left(
		(1-\cos d_{ij})\frac{\partial\theta_i^{jk}}{\partial u_j}
		+
		(1-\cos d_{ik})\frac{\partial\theta_i^{jk}}{\partial u_k}
		\right)=\sum_{j\sim i}B_{ij}(1-\cos d_{ij}).
	\end{aligned}
	\end{equation}
	Consequently, by \eqref{eq: derivative of K w.r.t. u} and \eqref{eq: Ai in Bij}, we have
	\[
	\frac{\partial K_i}{\partial u_i}
	=
	\sum_{j\sim i}B_{ij}
	-\sum_{j\sim i}B_{ij}(1-\cos d_{ij})
	=
	\sum_{j\sim i}B_{ij}\cos d_{ij}.
	\]
	The explicit expressions in \eqref{eq:spherical-curvature-A-B}
	and \eqref{eq:spherical-curvature-diagonal} follow by substituting
	$B_{ij}$ from \eqref{eq:spherical-curvature-off-diagonal} and using the identity
	$1-\cos d_{ij}=2\sin^2(d_{ij}/2)$.
\end{proof}

\subsection{Estimates for
	\texorpdfstring{$\epsilon$}{epsilon}-uniformly nondegenerate metrics}
We establish uniform estimates under the
$\epsilon$-uniform nondegeneracy condition and then control the change of interior angles under the spherical vertex scalings.
These results will be used in the proofs of short-time
existence and uniqueness for both combinatorial curvature flows.

For $0<\epsilon\leq\pi/2$, set
\[
C_1(\epsilon):=
\frac{\cot(\epsilon/2)}{\sin^2(\epsilon/2)}.
\]

\begin{lemma}\label{lem:spherical-uniform-bounds}
	Let $d=u*d_0$ be an $\epsilon$-uniformly nondegenerate
	PS metric on $\mathcal T$, and suppose that the vertex
	degrees are bounded by $D$.
	Then $|K_i|\leq(D+2)\pi$ for every $i\in V$. For every $i\in V$ and $j\sim i$, we have
\begin{equation}\label{eq:spherical-uniform-curvature-bounds}
	|B_{ij}|\leq C_1(\epsilon),\qquad
	|A_i|\leq2DC_1(\epsilon),\qquad
	\sum_{j\in V}
	\left|\frac{\partial K_i}{\partial u_j}\right|
	\leq2DC_1(\epsilon).
\end{equation}
\end{lemma}

\begin{proof}
	
	By the definition of the combinatorial curvature, we have
	\[
	|K_i|
	\leq2\pi+\sum_{\triangle ijk\in F}\theta_i^{jk}
	\leq(D+2)\pi.
	\]
	For each face $\triangle ijk\in F$, the definition of
	$\mathcal K_\epsilon^S$ gives
	\[
	-\pi+\epsilon
	\leq\theta_i^{jk}+\theta_j^{ik}-\theta_k^{ij}
	\leq\pi-\epsilon,
	\]
	and the same bounds hold after permuting the vertices.
	By \eqref{eq:spherical-cos-half-angle}, we have
	\[
	\cos^2\frac{d_{ij}}2
	=
	\frac{
		\cos\frac{\theta_i^{jk}+\theta_k^{ij}-\theta_j^{ik}}2
		\cos\frac{\theta_j^{ik}+\theta_k^{ij}-\theta_i^{jk}}2
	}{
		\sin\theta_i^{jk}\sin\theta_j^{ik}
	}.
	\]
	Both cosine factors are at least $\sin(\epsilon/2)$, and
	$0<\sin\theta_i^{jk}\sin\theta_j^{ik}\leq1$.
	Together with the angle bounds above, this gives
	\begin{equation}\label{eq:spherical-uniform-geometric-bounds}
		\left|
		\tan\frac{\theta_i^{jk}+\theta_j^{ik}-\theta_k^{ij}}2
		\right|
		\leq\cot\frac{\epsilon}{2},
		\qquad
		\cos^2\frac{d_{ij}}2
		\geq\sin^2\frac{\epsilon}{2}.
	\end{equation}
		By Lemma~\ref{lem:spherical-angle-variation},
\begin{equation}\label{eq:spherical-angle-derivative-bound}
    \left|\frac{\partial\theta_i^{jk}}{\partial u_j}\right|
    =
    \left|
    \frac{1}{2\cos^2(d_{ij}/2)}
    \tan\frac{\theta_i^{jk}+\theta_j^{ik}-\theta_k^{ij}}2
    \right|
    \leq
    \frac{\cot(\epsilon/2)}{2\sin^2(\epsilon/2)}
    =\frac{C_1(\epsilon)}2.
\end{equation}
By \eqref{eq:spherical-angle-derivative-bound} and the definition
of $B_{ij}$ in \eqref{eq:spherical-coefficients-definition},
we obtain
\[
\begin{aligned}
    |B_{ij}|
    &\leq
    \left|\frac{\partial\theta_i^{jk_1}}{\partial u_j}\right|
    +
    \left|\frac{\partial\theta_i^{jk_2}}{\partial u_j}\right|\leq
    \frac{\cot(\epsilon/2)}{\sin^2(\epsilon/2)}
    =C_1(\epsilon).
\end{aligned}
\]
	Proposition~\ref{prop:spherical-curvature-variation}
	then gives
	\[
	|A_i|
	\leq\sum_{j\sim i}|B_{ij}|(1-\cos d_{ij})
	\leq2DC_1(\epsilon),
	\]
	and
	\[
	\begin{aligned}
		\sum_{j\in V}
		\left|\frac{\partial K_i}{\partial u_j}\right|
		&=
		\left|\sum_{j\sim i}B_{ij}\cos d_{ij}\right|
		+\sum_{j\sim i}|B_{ij}|\leq2\sum_{j\sim i}|B_{ij}|
		\leq2DC_1(\epsilon).
	\end{aligned}
	\]
\end{proof}

The following lemma gives a uniform estimate for the change of interior angles under the spherical vertex scalings.

\begin{lemma}\label{lem:spherical_stability}
	Let $d$ be a PS metric on $\triangle ijk$ with interior
	angles $(\theta_i,\theta_j,\theta_k)\in\mathcal K_\epsilon^S$.
	There exists $\delta=\delta(\epsilon)>0$ with the following
	property. Suppose that
	\[
	\max\{|u_i|,|u_j|,|u_k|\}\leq\delta.
	\]
	Then $u*d$ is a PS metric on $\triangle ijk$.
	Its interior angles satisfy
	\begin{equation}\label{eq:angle_stability}
		|\widetilde\theta_r-\theta_r|
		\leq\frac{\epsilon}{12},
		\qquad r\in\{i,j,k\}.
	\end{equation}
\end{lemma}
The proof of Lemma \ref{lem:spherical_stability} is technical and not central to the main argument, so we put the proof in the appendix.

\begin{remark}\label{rem:spherical_stability}
	Let $d$ be an $\epsilon$-uniformly nondegenerate PS metric.
	If $\|u\|_{L^\infty(V)}\leq\delta(\epsilon)$, then $u*d$ is
	$\epsilon/2$-uniformly nondegenerate.
 For each face, \eqref{eq:angle_stability}
	gives
	\[
	\begin{aligned}
		\frac{\epsilon}{2}
		&\leq\widetilde\theta_r
		\leq\pi-\frac{\epsilon}{2},
		\qquad r\in\{i,j,k\},\\
		\widetilde\theta_i+\widetilde\theta_j
		+\widetilde\theta_k
		&\geq\pi+\epsilon-\frac{3\epsilon}{12}
		\geq\pi+\frac{\epsilon}{2},\\
		\widetilde\theta_i+\widetilde\theta_j
		-\widetilde\theta_k
		&\leq\pi-\epsilon+\frac{3\epsilon}{12}
		\leq\pi-\frac{\epsilon}{2}.
	\end{aligned}
	\]
	The last inequality also holds after permuting the vertices.
	Thus
	$(\widetilde\theta_i,\widetilde\theta_j,\widetilde\theta_k)
	\in\mathcal K_{\epsilon/2}^S$.
\end{remark}

\section{The spherical combinatorial Yamabe flow}

\subsection{Existence of the spherical combinatorial Yamabe flow}

We construct a solution of the combinatorial Yamabe flow on infinitely triangulated surfaces by solving finite-dimensional systems on finite subcomplexes and passing to a limit.
Let $(S,\mathcal T)$ be a connected, locally finite, infinitely
triangulated surface with an initial PS metric $d_0$.
The combinatorial distance $d_{\mathcal T}(i,j)$ between
$i,j\in V$ is the minimum number of edges in an edge path
joining $i$ to $j$.
Fix $i_0\in V$. For each integer $n\geq1$, set
\[
B_n(i_0)
=
\{j\in V:d_{\mathcal T}(i_0,j)\leq n\},
\]
and let $\mathcal T_n$ be the subcomplex consisting of all
simplices whose vertices lie in $B_n(i_0)$.
Denote its vertex, edge, and face sets by $V_n,E_n,F_n$,
respectively. In particular, $V_n=B_n(i_0)$.
Set
\[
\operatorname{int}(V_n)
=
\{i\in V_n:j\in V_n\text{ for every }j\sim i\},
\qquad
\partial V_n
=
V_n\setminus\operatorname{int}(V_n).
\]
By local finiteness and connectedness,
\[
|V_n|<\infty,\qquad
V_n\subseteq V_{n+1},\qquad
\bigcup_{n=1}^{\infty}V_n=V.
\]
Moreover, each fixed vertex $i$ belongs to
$\operatorname{int}(V_n)$ whenever
$n\geq d_{\mathcal T}(i_0,i)+1$.

We consider the following finite-dimensional system
\begin{equation}\label{eq:spherical-finite-flow}
	\begin{cases}
		\displaystyle
		\frac{du_i^{[n]}}{dt}=-K_i(u^{[n]}),
		& i\in\operatorname{int}(V_n),\\[6pt]
		u_i^{[n]}(t)=0,
		& i\in\partial V_n,\\[4pt]
		u_i^{[n]}(0)=0,
		& i\in V_n.
	\end{cases}
\end{equation}
Here $u^{[n]}$ is extended by zero to $V\setminus V_n$,
and $K_i(u^{[n]})$ is the combinatorial curvature of
$u^{[n]}*d_0$ on $\mathcal T$.
Only faces incident to vertices in $\operatorname{int}(V_n)$
are affected by these conformal factors, and there are
finitely many such faces.
Thus the conformal factors whose zero extensions define
PS metrics form an open neighborhood of the origin in
$\mathbb R^{\operatorname{int}(V_n)}$.
The right-hand side of \eqref{eq:spherical-finite-flow}
is smooth on this set.
Standard ODE theory therefore gives a local smooth solution.

\begin{proof}[Proof of Theorem~\ref{thm:vxste}(a)]
	Let $\delta_0=\delta(\epsilon)>0$ be given by
	Lemma~\ref{lem:spherical_stability}, and set
	\begin{equation}\label{eq:spherical-common-time}
		T_1:=\frac{\delta_0}{(D+2)\pi}.
	\end{equation}
	By Remark~\ref{rem:spherical_stability}, $u*d_0$ is an
	$\epsilon/2$-uniformly nondegenerate PS metric whenever
	$\|u\|_{L^\infty(V)}\leq\delta_0$.
	
	For each $n$, let $[0,\tau_n)$ be the maximal interval
	of existence of \eqref{eq:spherical-finite-flow}.
	Since $|K_i|\leq(D+2)\pi$ for every PS metric,
	the flow equation and the zero extension give
	\begin{equation}\label{eq:spherical-first-time-bound}
		\left|\frac{du_i^{[n]}}{dt}\right|
		\leq(D+2)\pi,
		\qquad i\in V,\quad 0\leq t<\tau_n.
	\end{equation}
	Integrating and using the zero initial condition, we obtain
	\begin{equation}\label{eq:spherical-u-bound}
		\|u^{[n]}(t)\|_{L^\infty(V)}
		\leq(D+2)\pi t\leq\delta_0,
		\qquad
		0\leq t<\min\{T_1,\tau_n\}.
	\end{equation}
	Suppose that $\tau_n\leq T_1$.
	By \eqref{eq:spherical-u-bound},
	\[
	\limsup_{t\uparrow\tau_n}
	\|u^{[n]}(t)\|_{L^\infty(V)}
	\leq(D+2)\pi\tau_n\leq\delta_0.
	\]
	The bound \eqref{eq:spherical-first-time-bound} implies
	that $u^{[n]}(t)$ has a limit as $t\uparrow\tau_n$.
	By Remark~\ref{rem:spherical_stability}, this limit
	defines an $\epsilon/2$-uniformly nondegenerate PS metric.
	We can therefore extend the solution beyond $\tau_n$
	by standard ODE theory, contradicting the maximality
	of $\tau_n$.
	Hence $\tau_n>T_1$ for every $n$.
	In particular, all finite-dimensional solutions are defined
	on $[0,T_1]$, and
	\eqref{eq:spherical-first-time-bound} and
	\eqref{eq:spherical-u-bound} hold throughout this interval.
	Their associated metrics are
	$\epsilon/2$-uniformly nondegenerate.
	
	For $i\in\operatorname{int}(V_n)$, the chain rule,
	\eqref{eq:spherical-first-time-bound}, and
	Lemma~\ref{lem:spherical-uniform-bounds}, with
	$\epsilon$ replaced by $\epsilon/2$, give
	\begin{equation}\label{eq:spherical-second-time-bound}
		\begin{aligned}
			\left|\frac{d^2u_i^{[n]}}{dt^2}\right|
			&=
			\left|
			\sum_{j\in V}
			\frac{\partial K_i}{\partial u_j}(u^{[n]})
			\frac{du_j^{[n]}}{dt}
			\right|\\
			&\leq
			(D+2)\pi
			\sum_{j\in V}
			\left|
			\frac{\partial K_i}{\partial u_j}(u^{[n]})
			\right|\\
			&\leq
			2D(D+2)\pi C_1(\epsilon/2).
		\end{aligned}
	\end{equation}
	For $i\notin\operatorname{int}(V_n)$, the second
	derivative is zero, so this estimate holds for all
	$i\in V$.
Combining \eqref{eq:spherical-first-time-bound},
\eqref{eq:spherical-u-bound}, and
\eqref{eq:spherical-second-time-bound}, we obtain
	\begin{equation}\label{eq:spherical-c2-bound}
		\|u_i^{[n]}\|_{C^2([0,T_1])}
		\leq
		\delta_0+(D+2)\pi
		+2D(D+2)\pi C_1(\epsilon/2).
	\end{equation}
	The right-hand side depends only on $\epsilon$ and $D$.

	Since $V$ is countable,  Arzel\`a--Ascoli theorem
	and a diagonal argument yield a subsequence
	$\{n_l\}_{l=1}^{\infty}$ such that
	\[
	u_i^{[n_l]}\longrightarrow u_i^*
	\qquad\text{in }C^1([0,T_1]),  \forall i\in V
	\]
	for some $u^*=(u^*_i)_{i\in V}$ with $u^*_i\in C^1([0,T_1])$ for every $i\in V$.
	Passing to the limit in \eqref{eq:spherical-u-bound}
	at each vertex gives
	\[
	\|u^*(t)\|_{L^\infty(V)}
	\leq(D+2)\pi t\leq\delta_0,
	\qquad t\in[0,T_1].
	\]
	By Remark~\ref{rem:spherical_stability},
	$u^*(t)*d_0$ is an $\epsilon/2$-uniformly nondegenerate
	PS metric.
	
	Fix $i\in V$.
	For all sufficiently large $l$,
	$i\in\operatorname{int}(V_{n_l})$.
	Since $K_i$ depends smoothly on the conformal factors
	at $i$ and its finitely many adjacent vertices,
	the above convergence and uniform nondegeneracy imply
	that $K_i(u^{[n_l]})$ converges to $K_i(u^*)$
	uniformly on $[0,T_1]$.
	Passing to the limit in \eqref{eq:spherical-finite-flow},
	we obtain
	\[
	\begin{cases}
		\dfrac{du_i^*}{dt}=-K_i(u^*),\\[4pt]
		u_i^*(0)=0,
	\end{cases}
	\qquad i\in V.
	\]
Since the right-hand side at each vertex depends smoothly on the conformal factors at only finitely many vertices, induction on the order of the time derivative gives \(u^*\in C_t^\infty(V\times[0,T])\).
\end{proof}

\subsection{Uniqueness of the spherical combinatorial Yamabe flow}
We prove the short-time uniqueness using
Lemma~\ref{lem:spherical-uniform-bounds}
and Gr\"onwall's inequality.

\begin{proof}[Proof of Theorem~\ref{thm:uq}(a)]
	Let $\delta_0$ and $T_1$ be as in the proof of Theorem \ref{thm:vxste}(a).
	The definition of combinatorial curvature gives
	$|K_i|\leq(D+2)\pi$ for every PS metric.
	Integrating the flow equations and using $u(0)=v(0)=0$, we obtain
	\[
	\max\left\{
	\|u(t)\|_{L^\infty(V)},\,
	\|v(t)\|_{L^\infty(V)}
	\right\}
	\leq(D+2)\pi t\leq\delta_0,
	\qquad t\in[0,T].
	\]
	Moreover, for $r,t\in[0,T]$,
	\[
	\|u(t)-u(r)\|_{L^\infty(V)}
	+\|v(t)-v(r)\|_{L^\infty(V)}
	\leq2(D+2)\pi|t-r|.
	\]
	It follows that
	$t\mapsto\|u(t)-v(t)\|_{L^\infty(V)}$
	is continuous on $[0,T]$.
	For $s\in[0,1]$, set
	\[
	w_s(t):=(1-s)v(t)+su(t).
	\]
	Then
	\[
	\begin{aligned}
		\|w_s(t)\|_{L^\infty(V)}
		&\leq(1-s)\|v(t)\|_{L^\infty(V)}
		+s\|u(t)\|_{L^\infty(V)}\leq\delta_0.
	\end{aligned}
	\]
	By 
	Remark~\ref{rem:spherical_stability}, the metric $w_s(t)*d_0$
	is $\epsilon/2$-uniformly nondegenerate for every
	$s\in[0,1]$ and $t\in[0,T]$.
	Applying Lemma~\ref{lem:spherical-uniform-bounds} with
	$\epsilon$ replaced by $\epsilon/2$, we obtain
	\begin{equation}\label{eq:uniqueness-row-bound}
		\sum_{j\in V}
		\left|\frac{\partial K_i}{\partial u_j}(w_s(t))\right|
		\leq2DC_1(\epsilon/2).
	\end{equation}
	Since $K_i$ depends only on the conformal factors at $i$
	and its finitely many adjacent vertices, the fundamental
	theorem of calculus gives
	\[
	\begin{aligned}
		K_i(u(t))-K_i(v(t))
		&=\int_0^1\frac{d}{ds}K_i(w_s(t))\,ds=\int_0^1\sum_{j\in V}
		\frac{\partial K_i}{\partial u_j}(w_s(t))
		\bigl(u_j(t)-v_j(t)\bigr)\,ds.
	\end{aligned}
	\]
	Using \eqref{eq:uniqueness-row-bound} and taking the supremum
	over $i\in V$, we obtain
	\begin{equation}\label{eq:uniqueness-curvature-lipschitz}
		\|K(u(t))-K(v(t))\|_{L^\infty(V)}
		\leq2DC_1(\epsilon/2)\|u(t)-v(t)\|_{L^\infty(V)}.
	\end{equation}
	The flow equations and the common initial value imply
	\[
	u_i(t)-v_i(t)
	=-\int_0^t\bigl(K_i(u(\tau))-K_i(v(\tau))\bigr)\,d\tau.
	\]
	Thus, by \eqref{eq:uniqueness-curvature-lipschitz},
	\[
	|u_i(t)-v_i(t)|
	\leq2DC_1(\epsilon/2)
	\int_0^t\|u(\tau)-v(\tau)\|_{L^\infty(V)}\,d\tau.
	\]
	Taking the supremum over $i\in V$ yields
	\[
	\|u(t)-v(t)\|_{L^\infty(V)}
	\leq2DC_1(\epsilon/2)
	\int_0^t\|u(\tau)-v(\tau)\|_{L^\infty(V)}\,d\tau.
	\]
Gr\"onwall's inequality therefore gives $u(t)=v(t)$
for all $t\in[0,T]$.
\end{proof}

\section{The spherical combinatorial Calabi flow}

\subsection{Existence of the spherical combinatorial Calabi flow}

We establish uniform bounds for the first- and second-order 
derivatives of the conformal factors and use these estimates
to prove short-time existence by an exhaustion argument.

\begin{lemma}\label{lem:spherical-calabi-first-derivative}
	Suppose that the PS metric $d=u*d_0$ is
	$\epsilon$-uniformly nondegenerate and that the vertex degrees
	of $\mathcal T$ are bounded by $D$.
	Then, along the spherical combinatorial Calabi flow,
	\begin{equation}\label{eq:spherical-calabi-first-derivative}
		\left|\frac{du_i}{dt}\right|
		\leq2D(D+2)\pi C_1(\epsilon),
		\qquad i\in V.
	\end{equation}
\end{lemma}

\begin{proof}
	By the Calabi flow equation and
	Lemma~\ref{lem:spherical-uniform-bounds},
	\[
	\left|\frac{du_i}{dt}\right|
	\leq
	\sum_{j\in V}
	\left|\frac{\partial K_i}{\partial u_j}\right||K_j|
	\leq
	(D+2)\pi
	\sum_{j\in V}
	\left|\frac{\partial K_i}{\partial u_j}\right|
	\leq2D(D+2)\pi C_1(\epsilon).
	\]
\end{proof}

To estimate the time derivatives of the coefficients, set
\[
\begin{aligned}
	C_B(\epsilon)
	&:=C_1(\epsilon)
	\left(
	2\cot^2\frac{\epsilon}{2}
	+\frac{3}{\sin^4(\epsilon/2)}
	\right),\\
	C_2(\epsilon,D)
	&:=2DC_B(\epsilon)+4DC_1(\epsilon).
\end{aligned}
\]

\begin{lemma}\label{lem:spherical-uniform-coefficient-variations}
	Suppose that the PS metric $d=u*d_0$ is
	$\epsilon$-uniformly nondegenerate and that the vertex degrees
	of $\mathcal T$ are bounded by $D$.
	Then, along the spherical combinatorial Calabi flow,
	\begin{equation}\label{eq:spherical-uniform-coefficient-variations}
		\begin{aligned}
			\left|\frac{dB_{ij}}{dt}\right|
			&\leq2D(D+2)\pi C_1(\epsilon)C_B(\epsilon),
			\qquad j\sim i,\\
			\left|\frac{dA_i}{dt}\right|
			&\leq2D(D+2)\pi C_1(\epsilon)C_2(\epsilon,D),\\
			\sum_{j\in V}
			\left|
			\frac{d}{dt}
			\left(\frac{\partial K_i}{\partial u_j}\right)
			\right|
			&\leq2D(D+2)\pi C_1(\epsilon)C_2(\epsilon,D).
		\end{aligned}
	\end{equation}
\end{lemma}
\begin{proof}
		For later use in the existence and uniqueness proofs,
		we first derive estimates for a general family of discrete
		conformal factors $u(t)$.
		Suppose that each $u_i$ is $C^1$ in $t$,
		$u_t(t)\in L^\infty(V)$, and $u(t)*d_0$ is
		$\epsilon$-uniformly nondegenerate.
		Differentiating \eqref{eq: vertex_scaling} gives
	\[
	\frac{d\,d_{ij}}{dt}
	=\tan\frac{d_{ij}}2
	\left(\frac{du_i}{dt}+\frac{du_j}{dt}\right).
	\]
	Hence
	\[
	\begin{aligned}
		\frac{d}{dt}
		\left(\frac{1}{2\cos^2(d_{ij}/2)}\right)
		&=
		\frac{\tan^2(d_{ij}/2)}{2\cos^2(d_{ij}/2)}
		\left(\frac{du_i}{dt}+\frac{du_j}{dt}\right),\\
		\frac{d}{dt}\cos d_{ij}
		&=
		-(1-\cos d_{ij})
		\left(\frac{du_i}{dt}+\frac{du_j}{dt}\right).
	\end{aligned}
	\]
Recall from \eqref{eq:spherical-curvature-off-diagonal} that
	\[
	B_{ij}
	=
	\frac{1}{2\cos^2(d_{ij}/2)}
	\sum_{\ell=1}^2
	\tan\frac{
		\theta_i^{jk_\ell}+\theta_j^{ik_\ell}
		-\theta_{k_\ell}^{ij}
	}{2}.
	\]
Applying the product and chain rules to this expression,
we obtain
	\begin{equation}\label{eq:spherical-B-time-derivative}
		\begin{aligned}
			\frac{dB_{ij}}{dt}
			&=
			\frac{d}{dt}
			\left(\frac{1}{2\cos^2(d_{ij}/2)}\right)
			\sum_{\ell=1}^2
			\tan\frac{
				\theta_i^{jk_\ell}+\theta_j^{ik_\ell}
				-\theta_{k_\ell}^{ij}
			}{2}\\
			&\quad+
			\frac{1}{2\cos^2(d_{ij}/2)}
			\sum_{\ell=1}^2
			\frac{d}{dt}
			\left(
			\tan\frac{
				\theta_i^{jk_\ell}+\theta_j^{ik_\ell}
				-\theta_{k_\ell}^{ij}
			}{2}
			\right)\\
			&=
			B_{ij}\tan^2\frac{d_{ij}}2
			\left(\frac{du_i}{dt}+\frac{du_j}{dt}\right)\\
			&\quad+
			\frac{1}{4\cos^2(d_{ij}/2)}
			\sum_{\ell=1}^2
			\frac{
				\dfrac{d\theta_i^{jk_\ell}}{dt}
				+\dfrac{d\theta_j^{ik_\ell}}{dt}
				-\dfrac{d\theta_{k_\ell}^{ij}}{dt}
			}{
				\cos^2\frac{
					\theta_i^{jk_\ell}+\theta_j^{ik_\ell}
					-\theta_{k_\ell}^{ij}
				}{2}
			}.
		\end{aligned}
	\end{equation}
Similarly, differentiating the identity
$A_i=\sum_{j\sim i}B_{ij}(1-\cos d_{ij})$
from Proposition~\ref{prop:spherical-curvature-variation} gives
	\begin{equation}\label{eq:spherical-A-time-derivative}
		\begin{aligned}
			\frac{dA_i}{dt}
			&=
			\sum_{j\sim i}
			(1-\cos d_{ij})\frac{dB_{ij}}{dt}+
			\sum_{j\sim i}B_{ij}(1-\cos d_{ij})
			\left(\frac{du_i}{dt}+\frac{du_j}{dt}\right).
		\end{aligned}
	\end{equation}
	By \eqref{eq:spherical-angle-derivative-bound},
\begin{equation}\label{eq:spherical-angle-partial-bounds-2}
    \left|\frac{\partial\theta_i^{jk}}{\partial u_j}\right|
    \leq\frac{C_1(\epsilon)}2,
    \qquad
    \left|\frac{\partial\theta_i^{jk}}{\partial u_k}\right|
    \leq\frac{C_1(\epsilon)}2.
\end{equation}
	The identity
	$\operatorname{Area}(\triangle ijk)
	=\theta_i^{jk}+\theta_j^{ik}+\theta_k^{ij}-\pi$,
	together with \eqref{eq:ddarea} and the symmetry of the
	angle derivatives, gives
	\begin{equation}\label{eq:spherical-angle-partial-ii}
	    \frac{\partial\theta_i^{jk}}{\partial u_i}
	=
	-\cos d_{ij}\frac{\partial\theta_i^{jk}}{\partial u_j}
	-\cos d_{ik}\frac{\partial\theta_i^{jk}}{\partial u_k}.
	\end{equation}
	Thus, by the chain rule, \eqref{eq:spherical-angle-partial-bounds-2} and \eqref{eq:spherical-angle-partial-ii}, we have
\begin{equation}\label{eq:spherical-angle-time-bound}
	\begin{aligned}
		\left|\frac{d\theta_i^{jk}}{dt}\right|
		&\leq
		\sum_{r\in\{i,j,k\}}
		\left|\frac{\partial\theta_i^{jk}}{\partial u_r}\right|
		\|u_t(t)\|_{L^\infty(V)}\\
		&\leq
		2\left(
		\left|\frac{\partial\theta_i^{jk}}{\partial u_j}\right|
		+\left|\frac{\partial\theta_i^{jk}}{\partial u_k}\right|
		\right)\|u_t(t)\|_{L^\infty(V)}\\
		&\leq2C_1(\epsilon)\|u_t(t)\|_{L^\infty(V)}.
	\end{aligned}
\end{equation}
	The same estimate holds for the other two angles.
	By \eqref{eq:spherical-uniform-geometric-bounds} and the $\epsilon$-uniformly nondegenerate condition,
    
	\begin{equation}\label{eq:spherical-trigonometric-bounds}
    \tan^2\frac{d_{ij}}2
    \leq\cot^2\frac{\epsilon}{2},
    \qquad
    \frac{1}{\cos^2(d_{ij}/2)}
    \leq\frac{1}{\sin^2(\epsilon/2)},
    \qquad
    \frac{1}{
        \cos^2\frac{
            \theta_i^{jk}+\theta_j^{ik}-\theta_k^{ij}
        }2
    }
    \leq\frac{1}{\sin^2(\epsilon/2)}.
\end{equation}
Using these bounds \eqref{eq:spherical-angle-time-bound}, \eqref{eq:spherical-trigonometric-bounds}, 
and the estimate $|B_{ij}|\leq C_1(\epsilon)$ established
in Lemma~\ref{lem:spherical-uniform-bounds}, we obtain
from \eqref{eq:spherical-B-time-derivative} that

\begin{equation}\label{eq:Bij_derivative-bound}
\begin{aligned}
\left|\frac{dB_{ij}}{dt}\right|
&\leq
2C_1(\epsilon)\cot^2\frac{\epsilon}{2}
\|u_t(t)\|_{L^\infty(V)}\\
&\quad+
\frac{1}{4\sin^4(\epsilon/2)}
\sum_{\ell=1}^2
\left(
\left|\frac{d\theta_i^{jk_\ell}}{dt}\right|
+\left|\frac{d\theta_j^{ik_\ell}}{dt}\right|
+\left|\frac{d\theta_{k_\ell}^{ij}}{dt}\right|
\right)\\
&\leq
C_1(\epsilon)
\left(
2\cot^2\frac{\epsilon}{2}
+\frac{3}{\sin^4(\epsilon/2)}
\right)\|u_t(t)\|_{L^\infty(V)}\\
&=C_B(\epsilon)\|u_t(t)\|_{L^\infty(V)}.
\end{aligned}
\end{equation}
Since $0<1-\cos d_{ij}\leq2$ and $\deg(i)\leq D$,
	\eqref{eq:spherical-A-time-derivative} gives
    
\begin{equation}\label{eq:Ai_derivative-bound}
\begin{aligned}
\left|\frac{dA_i}{dt}\right|
&\leq
2\sum_{j\sim i}\left|\frac{dB_{ij}}{dt}\right|
+2\sum_{j\sim i}|B_{ij}|
\left(
\left|\frac{du_i}{dt}\right|
+\left|\frac{du_j}{dt}\right|
\right)\\
&\leq
\bigl(2DC_B(\epsilon)+4DC_1(\epsilon)\bigr)
\|u_t(t)\|_{L^\infty(V)}\\
&=C_2(\epsilon,D)\|u_t(t)\|_{L^\infty(V)}.
\end{aligned}
\end{equation}
Since $\partial K_i/\partial u_j=0$ for $j\ne i$ and
$j\not\sim i$, Proposition~\ref{prop:spherical-curvature-variation}
gives
\[
\begin{aligned}
	\sum_{j\in V}
	\left|
	\frac{d}{dt}
	\left(\frac{\partial K_i}{\partial u_j}\right)
	\right|
	&=
	\left|
	\frac{d}{dt}
	\left(\frac{\partial K_i}{\partial u_i}\right)
	\right|
	+
	\sum_{j\sim i}
	\left|
	\frac{d}{dt}
	\left(\frac{\partial K_i}{\partial u_j}\right)
	\right|\\
	&=
	\left|
	\sum_{j\sim i}
	\left(
	\cos d_{ij}\frac{dB_{ij}}{dt}
	+B_{ij}\frac{d}{dt}\cos d_{ij}
	\right)
	\right|
	+\sum_{j\sim i}\left|\frac{dB_{ij}}{dt}\right|\\
	&\leq
	2\sum_{j\sim i}\left|\frac{dB_{ij}}{dt}\right|
	+\sum_{j\sim i}|B_{ij}|
	\left|\frac{d}{dt}\cos d_{ij}\right|.
\end{aligned}
\]
The last inequality follows from the triangle inequality
and $|\cos d_{ij}|\leq1$.
	Since
	\[
	\left|\frac{d}{dt}\cos d_{ij}\right|
	\leq
	(1-\cos d_{ij})
	\left(
	\left|\frac{du_i}{dt}\right|
	+\left|\frac{du_j}{dt}\right|
	\right)
	\leq4\|u_t(t)\|_{L^\infty(V)},
	\]
	we conclude that
\begin{equation}\label{eq:spherical-general-jacobian-bound}
    \begin{aligned}
        \sum_{j\in V}
        \left|
        \frac{d}{dt}
        \left(\frac{\partial K_i}{\partial u_j}(u(t))\right)
        \right|
        &\leq
        \bigl(2DC_B(\epsilon)+4DC_1(\epsilon)\bigr)
        \|u_t(t)\|_{L^\infty(V)}\\
        &=C_2(\epsilon,D)\|u_t(t)\|_{L^\infty(V)}.
    \end{aligned}
\end{equation}
Along the spherical combinatorial Calabi flow,
Lemma~\ref{lem:spherical-calabi-first-derivative} gives
\[
\|u_t(t)\|_{L^\infty(V)}
\leq2D(D+2)\pi C_1(\epsilon).
\]
Substituting this bound into
\eqref{eq:Bij_derivative-bound},
\eqref{eq:Ai_derivative-bound}, and
\eqref{eq:spherical-general-jacobian-bound} yields
\eqref{eq:spherical-uniform-coefficient-variations}.
\end{proof}

We next estimate the second-order derivatives along the spherical combinatorial Calabi flow.

\begin{proposition}\label{prop:spherical-uniform-flow-estimates}
	Suppose that the PS metric $d=u*d_0$ is
	$\epsilon$-uniformly nondegenerate and that the vertex degrees
	of $\mathcal T$ are bounded by $D$.
	Then, along the spherical combinatorial Calabi flow,
	\begin{equation}\label{eq:spherical-calabi-uniform-time-estimates}
		\begin{aligned}
			\left|\frac{d^2u_i}{dt^2}\right|
			&\leq2D(D+2)\pi C_1(\epsilon)
			\left(
			4D^2C_1(\epsilon)^2
			+(D+2)\pi C_2(\epsilon,D)
			\right),
			\qquad i\in V.
		\end{aligned}
	\end{equation}
\end{proposition}

\begin{proof}
	By the chain rule and
	Lemma~\ref{lem:spherical-uniform-bounds},
	\[
	\left|\frac{dK_i}{dt}\right|
	\leq
	\sum_{j\in V}
	\left|\frac{\partial K_i}{\partial u_j}\right|
	\left|\frac{du_j}{dt}\right|
	\leq2DC_1(\epsilon)\|u_t\|_{L^\infty(V)}.
	\]
	Differentiating the spherical combinatorial Calabi flow equation gives
	\[
	\frac{d^2u_i}{dt^2}
	=
	-\sum_{j\in V}
	\frac{d}{dt}
	\left(\frac{\partial K_i}{\partial u_j}\right)K_j
	-\sum_{j\in V}
	\frac{\partial K_i}{\partial u_j}\frac{dK_j}{dt}.
	\]
Using Lemma~\ref{lem:spherical-uniform-bounds} and
\eqref{eq:spherical-general-jacobian-bound}, we obtain
	\[
	\begin{aligned}
		\left|\frac{d^2u_i}{dt^2}\right|
		&\leq
		(D+2)\pi
		\sum_{j\in V}
		\left|
		\frac{d}{dt}
		\left(\frac{\partial K_i}{\partial u_j}\right)
		\right|
		+2DC_1(\epsilon)\|K_t\|_{L^\infty(V)}\\
		&\leq
\bigl((D+2)\pi C_2(\epsilon,D)+4D^2C_1(\epsilon)^2\bigr)
\|u_t\|_{L^\infty(V)}.
	\end{aligned}
	\]
	The conclusion follows from
	Lemma~\ref{lem:spherical-calabi-first-derivative}.
\end{proof}

We use the finite exhaustion $\mathcal T_n$ constructed
in Section~3.1 and extend $u^{[n]}$ by zero to
$V\setminus V_n$.
Consider the finite-dimensional system
\begin{equation}\label{eq:spherical-calabi-finite-system}
	\begin{cases}
		\displaystyle
		\frac{du_i^{[n]}}{dt}
		=-\sum_{j\in V}
		\frac{\partial K_i}{\partial u_j}(u^{[n]})
		K_j(u^{[n]}),
		& i\in\operatorname{int}(V_n),\\[6pt]
		u_i^{[n]}(t)=0,
		& i\in\partial V_n,\\[4pt]
		u_i^{[n]}(0)=0,
		& i\in V_n.
	\end{cases}
\end{equation}
Here the curvatures and their partial derivatives are
computed on $\mathcal T$ and evaluated at $u^{[n]}$.
As in Section~3.1, standard ODE theory gives a local
smooth solution for each $n$.
\begin{proof}[Proof of Theorem~\ref{thm:vxste}(b)]
	Let $\delta_0$ be as in the proof of
	Theorem~\ref{thm:vxste}(a), and set
	\begin{equation}\label{eq:spherical-calabi-existence-time}
		T_2:=
		\frac{\delta_0}{2D(D+2)\pi C_1(\epsilon/2)}.
	\end{equation}
	For each $n$, let $[0,\tau_n)$ be the maximal interval
	of existence of \eqref{eq:spherical-calabi-finite-system}.
	
	As long as $\|u^{[n]}(t)\|_{L^\infty(V)}\leq\delta_0$,
	Remark~\ref{rem:spherical_stability} implies that
	$u^{[n]}(t)*d_0$ is $\epsilon/2$-uniformly nondegenerate.
	The argument in
	Lemma~\ref{lem:spherical-calabi-first-derivative},
	with $\epsilon$ replaced by $\epsilon/2$, applies at
	every vertex in $\operatorname{int}(V_n)$.
	Since the time derivatives vanish at all remaining
	vertices, we obtain
	\begin{equation}\label{eq:spherical-calabi-finite-first-bound}
		\left|\frac{du_i^{[n]}}{dt}\right|
		\leq2D(D+2)\pi C_1(\epsilon/2).
	\end{equation}
	Integrating and using $u^{[n]}(0)=0$ gives
	\[
	\|u^{[n]}(t)\|_{L^\infty(V)}
	\leq2D(D+2)\pi C_1(\epsilon/2)t.
	\]
	The norm therefore cannot first reach $\delta_0$
	before $T_2$, so both estimates hold for
	$0\leq t<\min\{T_2,\tau_n\}$.
	If $\tau_n\leq T_2$, the derivative bound implies that
	$u^{[n]}(t)$ has a limit as $t\uparrow\tau_n$,
	with norm at most $\delta_0$.
	The same continuation argument as in the Yamabe case
	then extends the solution beyond $\tau_n$,
	contradicting maximality.
	Thus every $u^{[n]}$ is defined on $[0,T_2]$,
	and the above estimates hold throughout this interval.
	In particular, the associated metrics are
	$\epsilon/2$-uniformly nondegenerate.
	
	Applying \eqref{eq:spherical-general-jacobian-bound}
	to the zero extension $u^{[n]}$, with $\epsilon$
	replaced by $\epsilon/2$, and repeating the calculations
	in the proof of
	Proposition~\ref{prop:spherical-uniform-flow-estimates},
	we obtain
	\[
	\left|\frac{d^2u_i^{[n]}}{dt^2}\right|
	\leq
	\Bigl(
	4D^2C_1(\epsilon/2)^2
	+(D+2)\pi C_2(\epsilon/2,D)
	\Bigr)
	\|u_t^{[n]}\|_{L^\infty(V)}
	\]
	for $i\in\operatorname{int}(V_n)$.
	At all remaining vertices, the second derivative is zero.
	Together with the preceding bounds, this gives a bound
	for $\|u_i^{[n]}\|_{C^2([0,T_2])}$ depending only on
	$\epsilon$ and $D$, independently of $n$ and $i$.
	
	As in the proof of Theorem~\ref{thm:vxste}(a),
	the Arzel\`a--Ascoli theorem and a diagonal argument
	yield a subsequence $\{n_l\}$ such that
	\[
	u_i^{[n_l]}\longrightarrow u_i^*
	\qquad\text{in }C^1([0,T_2]), \forall i\in V
	\]
	for some $u^*=(u^*_i)_{i\in V}$ with $u^*_i\in C^1([0,T_2])$ for every $i\in V$..
	The limit satisfies $u^*(0)=0$ and
	$\|u^*(t)\|_{L^\infty(V)}\leq\delta_0$.
	Hence $u^*(t)*d_0$ is
	$\epsilon/2$-uniformly nondegenerate.
	
	For each fixed $i\in V$, we have
	$i\in\operatorname{int}(V_{n_l})$ for all sufficiently
	large $l$.
	The right-hand side of the equation at $i$ depends
	smoothly on the conformal factors at the finitely many
	vertices of combinatorial distance at most two from $i$.
Passing to the limit in
\eqref{eq:spherical-calabi-finite-system}, we obtain
\[
\begin{cases}
	\displaystyle
	\frac{du_i^*}{dt}
	=-\sum_{j\in V}
	\frac{\partial K_i}{\partial u_j}(u^*)K_j(u^*),\\[6pt]
	u_i^*(0)=0,
\end{cases}
\qquad i\in V.
\]
Thus $u^*$ satisfies \eqref{Eq:VCLBF}.
	As in the case of combinatorial Yamabe flow, repeated differentiation gives
	$u^*\in C_t^\infty(V\times[0,T_2])$.
\end{proof}

\begin{remark}\label{rem:spherical-calabi-continuity}
	By \eqref{eq:spherical-calabi-finite-first-bound}
	and the convergence at each vertex, the solution
	constructed above satisfies
	\[
	\|u^*(t)-u^*(s)\|_{L^\infty(V)}
	\leq2D(D+2)\pi C_1(\epsilon/2)|t-s|,
	\qquad s,t\in[0,T_2].
	\]
	Hence $u^*\in C([0,T_2];L^\infty(V))$,
	so it satisfies the continuity assumption in
	Theorem~\ref{thm:uq}(b).
\end{remark}

\subsection{Uniqueness of the spherical combinatorial Calabi flow}
We prove uniqueness using the estimates above and Gr\"onwall's inequality.

\begin{proof}[Proof of Theorem~\ref{thm:uq}(b)]
	Let $u$ and $v$ be the two solutions in the theorem, and let
	$\delta_0$ and $T_2$ be as in the proof of
	Theorem~\ref{thm:vxste}(b).
	We first show that
	\[
	\max\{\|u(t)\|_{L^\infty(V)},\|v(t)\|_{L^\infty(V)}\}
	\leq\delta_0,\qquad t\in[0,T].
	\]
	Suppose otherwise. Since $u$ and $v$ are continuous in the
	$L^\infty$ norm and $u(0)=v(0)=0$, there exist
	$0<t_0<t_1\leq T$ such that the maximum first equals
	$\delta_0$ at $t_0$ and exceeds $\delta_0$ at $t_1$.
	By Remark~\ref{rem:spherical_stability}, both metrics are
	$\epsilon/2$-uniformly nondegenerate on $[0,t_0]$.
	Integrating the estimate in
	Lemma~\ref{lem:spherical-calabi-first-derivative}, with
	$\epsilon$ replaced by $\epsilon/2$, gives
	\[
	\begin{aligned}
		\delta_0
		&=\max\{\|u(t_0)\|_{L^\infty(V)},
		\|v(t_0)\|_{L^\infty(V)}\}\leq2D(D+2)\pi C_1(\epsilon/2)t_0<\delta_0,
	\end{aligned}
	\]
	where $t_0<T\leq T_2$ is used in the last inequality.
	This is a contradiction.
	
	Fix $t\in[0,T]$ and set
	$w_s=(1-s)v(t)+su(t)$ for $s\in[0,1]$.
	Since $\|w_s\|_{L^\infty(V)}\leq\delta_0$,
	Remark~\ref{rem:spherical_stability} implies that $w_s*d_0$
	is $\epsilon/2$-uniformly nondegenerate.
	For a discrete conformal factor $w$, let $\Delta_w$ denote
	the discrete Laplace operator associated with $w*d_0$.
	Then
	\[
	(\Delta_{w_s}K(w_s))_i
	=-\sum_{j\in V}\frac{\partial K_i}{\partial u_j}(w_s)K_j(w_s).
	\]
	All sums below have only finitely many nonzero terms for
	each fixed $i$. Since
	\[
	\frac{dw_s}{ds}=u(t)-v(t),
	\]
	applying \eqref{eq:spherical-general-jacobian-bound} with
	$t$ replaced by $s$ and $\epsilon$ by $\epsilon/2$ gives
	\begin{equation}\label{eq:spherical-jacobian-interpolation-bound}
    \sum_{j\in V}
    \left|
    \frac{d}{ds}
    \left(\frac{\partial K_i}{\partial u_j}(w_s)\right)
    \right|
    \leq C_2(\epsilon/2,D)
    \|u(t)-v(t)\|_{L^\infty(V)}.
\end{equation}
	Hence, by the chain rule,
	Lemma~\ref{lem:spherical-uniform-bounds} and \eqref{eq:spherical-jacobian-interpolation-bound},
	\[
	\begin{aligned}
		&\left|\frac{d}{ds}(\Delta_{w_s}K(w_s))_i\right|\\
		&\quad\leq\sum_{j\in V}
		\left|\frac{d}{ds}
		\left(\frac{\partial K_i}{\partial u_j}(w_s)\right)\right|
		|K_j(w_s)|\\
		&\qquad+\sum_{j,k\in V}
		\left|\frac{\partial K_i}{\partial u_j}(w_s)\right|
		\left|\frac{\partial K_j}{\partial u_k}(w_s)\right|
		|u_k(t)-v_k(t)|\\
		&\quad\leq
		\left(4D^2C_1(\epsilon/2)^2
		+(D+2)\pi C_2(\epsilon/2,D)\right)
		\|u(t)-v(t)\|_{L^\infty(V)}.
	\end{aligned}
	\]
Integrating with respect to $s$ and taking the supremum
	over $i\in V$, we obtain
	\begin{equation}\label{eq:spherical-calabi-rhs-lipschitz}
		\begin{aligned}
			&\|\Delta_{u(t)}K(u(t))-\Delta_{v(t)}K(v(t))\|_{L^\infty(V)}\\
			&\quad\leq
			\left(4D^2C_1(\epsilon/2)^2
			+(D+2)\pi C_2(\epsilon/2,D)\right)
			\|u(t)-v(t)\|_{L^\infty(V)}.
		\end{aligned}
	\end{equation}
	The flow equations and the common initial value give
	\[
	u_i(t)-v_i(t)
	=\int_0^t
	\left((\Delta_{u(\tau)}K(u(\tau)))_i
	-(\Delta_{v(\tau)}K(v(\tau)))_i\right)\,d\tau.
	\]
	Using \eqref{eq:spherical-calabi-rhs-lipschitz} and taking
	the supremum over $i\in V$, we obtain
	\[
	\begin{aligned}
		&\|u(t)-v(t)\|_{L^\infty(V)}\\
		&\quad\leq
		\left(4D^2C_1(\epsilon/2)^2
		+(D+2)\pi C_2(\epsilon/2,D)\right)
		\int_0^t\|u(\tau)-v(\tau)\|_{L^\infty(V)}\,d\tau.
	\end{aligned}
	\]
	Since $t\mapsto\|u(t)-v(t)\|_{L^\infty(V)}$ is continuous,
	Gr\"onwall's inequality gives $u(t)=v(t)$ on $[0,T]$.
\end{proof}

\section{Appendix}

\begin{proof}[Proof of Lemma~\ref{lem:spherical_stability}]
	We first establish uniform bounds for the side lengths.
	The assumptions imply $0<\epsilon\leq\pi/2$.
	By the definition of $\mathcal K_\epsilon^S$, we have
	\[
	\pi+\epsilon
	\leq\theta_i+\theta_j+\theta_k
	\leq3\pi-\epsilon,
	\qquad
	-\pi+\epsilon
	\leq\theta_i+\theta_j-\theta_k
	\leq\pi-\epsilon,
	\]
	and the latter inequalities hold cyclically. Hence
	\[
	-\cos\frac{\theta_i+\theta_j+\theta_k}{2}
	\geq\sin\frac{\epsilon}{2},
	\qquad
	\cos\frac{\theta_i+\theta_j-\theta_k}{2}
	\geq\sin\frac{\epsilon}{2}.
	\]
	By \eqref{eq:spherical-half-angle}, its cyclic permutations, and \eqref{eq:spherical-cos-half-angle},
	\begin{align*}
		\sin^2\frac{d_{ij}}2
		=
		-\frac{
			\cos\frac{\theta_i+\theta_j+\theta_k}{2}
			\cos\frac{\theta_i+\theta_j-\theta_k}{2}
		}{\sin\theta_i\sin\theta_j}
		\geq\sin^2\frac{\epsilon}{2},
	\end{align*}
	and
	\begin{align*}
		\cos^2\frac{d_{ij}}2
	=
		\frac{
			\cos\frac{\theta_i+\theta_k-\theta_j}{2}
			\cos\frac{\theta_j+\theta_k-\theta_i}{2}
		}{\sin\theta_i\sin\theta_j}
		\geq\sin^2\frac{\epsilon}{2}.
	\end{align*}
	Here we used $0<\sin\theta_i\sin\theta_j\leq1$.
	The same argument applies to the other two sides. Thus
	\begin{equation}\label{eq:spherical_initial_side_bounds}
		\epsilon\leq d_{rs}\leq\pi-\epsilon,
		\qquad \{r,s\}\subset\{i,j,k\}.
	\end{equation}
	
	We next estimate the differences in the triangle inequalities
	and the perimeter condition. By the spherical cosine law,
	\[
	\cos d_{jk}
	=
	\cos d_{ij}\cos d_{ik}
	+\sin d_{ij}\sin d_{ik}\cos\theta_i.
	\]
	Hence
	\begin{align*}
		&\sin\frac{d_{ij}+d_{ik}+d_{jk}}2
		\sin\frac{d_{ij}+d_{ik}-d_{jk}}2\\
		&\qquad=
		\frac{\cos d_{jk}-\cos(d_{ij}+d_{ik})}{2}\\
		&\qquad=
		\frac{
			\cos d_{jk}-\cos d_{ij}\cos d_{ik}
			+\sin d_{ij}\sin d_{ik}
		}{2}\\
		&\qquad=
		\frac12\sin d_{ij}\sin d_{ik}(1+\cos\theta_i)\\
		&\qquad=
		\sin d_{ij}\sin d_{ik}\cos^2\frac{\theta_i}{2}\\
		&\qquad\geq
		\sin^2\epsilon\sin^2\frac{\epsilon}{2},
	\end{align*}
	where the last inequality follows from
	\eqref{eq:spherical_initial_side_bounds} and
	$\theta_i\leq\pi-\epsilon$.
	Both sine factors are positive and at most one, so each is
	bounded below by $\sin^2\epsilon\sin^2(\epsilon/2)$.
	Using $\sin x\leq x$ for $x\geq0$ and
	\[
	\sin\frac{d_{ij}+d_{ik}+d_{jk}}2
	=
	\sin\frac{2\pi-d_{ij}-d_{ik}-d_{jk}}2,
	\]
	we obtain
	\begin{equation}\label{eq:spherical_initial_margins}
		\begin{aligned}
			d_{ij}+d_{ik}-d_{jk}
			&\geq2\sin^2\epsilon\sin^2\frac{\epsilon}{2},\\
			2\pi-d_{ij}-d_{ik}-d_{jk}
			&\geq2\sin^2\epsilon\sin^2\frac{\epsilon}{2}.
		\end{aligned}
	\end{equation}
	The first inequality also holds after cyclic permutation.
	
	Choose $\delta>0$ sufficiently small such that
	\begin{equation}\label{eq:spherical_scaling_range}
		e^\delta\cos\frac{\epsilon}{2}
		\leq\cos\frac{\epsilon}{4}.
	\end{equation}
	If $\max\{|u_i|,|u_j|,|u_k|\}\leq\delta$, then
	\eqref{eq:spherical_initial_side_bounds} and
	\eqref{eq:spherical_scaling_range} give
	\[
	0<e^{-\delta}\sin\frac{\epsilon}{2}
	\leq e^{(u_r+u_s)/2}\sin\frac{d_{rs}}2
	\leq e^\delta\cos\frac{\epsilon}{2}
	\leq\cos\frac{\epsilon}{4}<1.
	\]
	Therefore, the scaled lengths are uniquely defined in $(0,\pi)$ by
	\[
	\widetilde d_{rs}
	=
	2\arcsin\left(
	e^{(u_r+u_s)/2}\sin\frac{d_{rs}}2
	\right).
	\]
	For $0\leq x\leq\cos(\epsilon/4)$, we have
	\[
	\frac{d}{dx}\arcsin x
	=\frac{1}{\sqrt{1-x^2}}
	\leq\frac{1}{\sin(\epsilon/4)}.
	\]
	Since both $\sin(d_{rs}/2)$ and
	$e^{(u_r+u_s)/2}\sin(d_{rs}/2)$ lie in
	$[0,\cos(\epsilon/4)]$, the mean value theorem gives
	\begin{equation}\label{eq:uniform_side_perturbation}
		\begin{aligned}
			|\widetilde d_{rs}-d_{rs}|
			&=
			2\left|
			\arcsin\left(
			e^{(u_r+u_s)/2}\sin\frac{d_{rs}}2
			\right)
			-
			\arcsin\left(\sin\frac{d_{rs}}2\right)
			\right|\\
			&\leq
			\frac{2}{\sin(\epsilon/4)}
			\left|e^{(u_r+u_s)/2}-1\right|
			\sin\frac{d_{rs}}2\\
			&\leq
			\frac{2\cos(\epsilon/2)}{\sin(\epsilon/4)}
			(e^\delta-1).
		\end{aligned}
	\end{equation}
	Decrease $\delta$ if necessary so that
	\begin{equation}\label{eq:spherical_side_smallness}
		\frac{2\cos(\epsilon/2)}{\sin(\epsilon/4)}
		(e^\delta-1)
		\leq
		\min\left\{
		\frac{\epsilon}{2},\,
		\frac13\sin^2\epsilon\sin^2\frac{\epsilon}{2}
		\right\}.
	\end{equation}
	It follows from \eqref{eq:spherical_initial_side_bounds},
	\eqref{eq:uniform_side_perturbation}, and
	\eqref{eq:spherical_side_smallness} that
	\begin{equation}\label{eq:spherical_scaled_side_bounds}
		\frac{\epsilon}{2}
		\leq\widetilde d_{rs}
		\leq\pi-\frac{\epsilon}{2}.
	\end{equation}
	Moreover, each expression on the left-hand side of
	\eqref{eq:spherical_initial_margins} changes by at most
	\[
	|\widetilde d_{ij}-d_{ij}|
	+
	|\widetilde d_{ik}-d_{ik}|
	+
	|\widetilde d_{jk}-d_{jk}|
	\leq\sin^2\epsilon\sin^2\frac{\epsilon}{2}.
	\]
	Consequently,
	\[
	\widetilde d_{ij}+\widetilde d_{ik}-\widetilde d_{jk}
	\geq\sin^2\epsilon\sin^2\frac{\epsilon}{2}>0
	\quad\text{(cyclically)},
	\]
	and
	\[
	2\pi-\widetilde d_{ij}-\widetilde d_{ik}-\widetilde d_{jk}
	\geq\sin^2\epsilon\sin^2\frac{\epsilon}{2}>0.
	\]
	Thus $u*d$ is a PS metric on $\triangle ijk$.
	
	We now estimate the variation of the interior angles.
By the spherical cosine law,
\[
\begin{aligned}
	\cos d_{jk}
	&=
	\cos d_{ij}\cos d_{ik}
	+\sin d_{ij}\sin d_{ik}\cos\theta_i,\\
	\cos\widetilde d_{jk}
	&=
	\cos\widetilde d_{ij}\cos\widetilde d_{ik}
	+\sin\widetilde d_{ij}\sin\widetilde d_{ik}
	\cos\widetilde\theta_i.
\end{aligned}
\]
Subtracting these identities and rearranging, we obtain
\begin{equation}\label{eq:spherical-cosine-difference-identity}
	\begin{aligned}
		&(\cos\widetilde\theta_i-\cos\theta_i)
		\sin\widetilde d_{ij}\sin\widetilde d_{ik}\\
		&\quad=
		\cos\widetilde d_{jk}-\cos d_{jk}
		-\bigl(
		\cos\widetilde d_{ij}\cos\widetilde d_{ik}
		-\cos d_{ij}\cos d_{ik}
		\bigr)\\
		&\qquad\quad-
		\cos\theta_i
		\bigl(
		\sin\widetilde d_{ij}\sin\widetilde d_{ik}
		-\sin d_{ij}\sin d_{ik}
		\bigr).
	\end{aligned}
\end{equation}
Using $|\cos x-\cos y|\leq|x-y|$ and $|\cos x|\leq1$,
	we obtain
\begin{equation}\label{eq:spherical-cosine-product-bound}
	\begin{aligned}
		&\left|
		\cos\widetilde d_{ij}\cos\widetilde d_{ik}
		-\cos d_{ij}\cos d_{ik}
		\right|\\
		&\quad\leq
		\left|\cos\widetilde d_{ij}-\cos d_{ij}\right|
		\left|\cos\widetilde d_{ik}\right|
		+
		\left|\cos d_{ij}\right|
		\left|\cos\widetilde d_{ik}-\cos d_{ik}\right|\\
		&\quad\leq
		|\widetilde d_{ij}-d_{ij}|
		+|\widetilde d_{ik}-d_{ik}|.
	\end{aligned}
\end{equation}
	Similarly, using $|\sin x-\sin y|\leq|x-y|$ and
	$|\sin x|\leq1$, we have
\begin{equation}\label{eq:spherical-sine-product-bound}
	\left|
	\sin\widetilde d_{ij}\sin\widetilde d_{ik}
	-\sin d_{ij}\sin d_{ik}
	\right|
	\leq
	|\widetilde d_{ij}-d_{ij}|
	+|\widetilde d_{ik}-d_{ik}|.
\end{equation}
	Moreover, \eqref{eq:spherical_scaled_side_bounds} gives
\begin{equation}\label{eq:spherical-sine-product-lower-bound}
	\sin\widetilde d_{ij}\sin\widetilde d_{ik}
	\geq\sin^2\frac{\epsilon}{2}.
\end{equation}
Taking absolute values in
\eqref{eq:spherical-cosine-difference-identity} and using
\eqref{eq:spherical-cosine-product-bound}, \eqref{eq:spherical-sine-product-bound},
\eqref{eq:spherical-sine-product-lower-bound}, and
$|\cos\theta_i|\leq1$, we obtain
\begin{equation}\label{eq:spherical_cosine_perturbation}
	\begin{aligned}
		|\cos\widetilde\theta_i-\cos\theta_i|
		&\leq
		\frac{
			|\widetilde d_{jk}-d_{jk}|
			+2|\widetilde d_{ij}-d_{ij}|
			+2|\widetilde d_{ik}-d_{ik}|
		}{\sin^2(\epsilon/2)}\\
		&\leq
		\frac{10\cos(\epsilon/2)}
		{\sin^2(\epsilon/2)\sin(\epsilon/4)}
		(e^\delta-1),
	\end{aligned}
\end{equation}
where the last inequality follows from
\eqref{eq:uniform_side_perturbation}.
The same estimate holds for the other two angles.

	Finally, decrease $\delta$ if necessary so that
	\[
	\frac{10\cos(\epsilon/2)}
	{\sin^2(\epsilon/2)\sin(\epsilon/4)}
	(e^\delta-1)
	\leq1-\cos\frac{\epsilon}{12}.
	\]
	By \eqref{eq:spherical_cosine_perturbation}, we have
	\[
	|\cos\widetilde\theta_r-\cos\theta_r|
	\leq1-\cos\frac{\epsilon}{12},
	\qquad r\in\{i,j,k\}.
	\]
	For any $a,b\in[0,\pi]$, we have
	\[
	\frac{a+b}{2}
	\in
	\left[
	\frac{|a-b|}{2},\,
	\pi-\frac{|a-b|}{2}
	\right].
	\]
	Therefore,
	\[
	\begin{aligned}
		|\cos a-\cos b|
		&=
		2\sin\frac{a+b}{2}\sin\frac{|a-b|}{2}\\
		&\geq
		2\sin^2\frac{|a-b|}{2}
		=
		1-\cos|a-b|.
	\end{aligned}
	\]
	Applying this inequality to
	$a=\widetilde\theta_r$ and $b=\theta_r$ gives
	\[
	1-\cos|\widetilde\theta_r-\theta_r|
	\leq1-\cos\frac{\epsilon}{12}.
	\]
	Since $1-\cos x$ is strictly increasing on $[0,\pi]$, we obtain
	\[
	|\widetilde\theta_r-\theta_r|
	\leq\frac{\epsilon}{12},
	\qquad r\in\{i,j,k\}.
	\]
All the above restrictions can be satisfied by choosing
$\delta=\delta(\epsilon)>0$ sufficiently small.
\end{proof}

\end{document}